\documentclass[11pt]{article}
 \usepackage{euscript,amsmath,amsfonts,amssymb,epsfig}
 \title{The example of a self-similar continuum which is not an attractor of any zipper.}
 \author{Purevdorj O., Tetenov A.V. }
 \date{}
  \begin{document}
\maketitle

 \newcommand {\al} {\alpha}
\newcommand {\be} {\beta}
\newcommand {\da} {\delta}
\newcommand {\Da} {\Delta}
\newcommand {\ga} {\gamma}
\newcommand {\Ga} {\Gamma}

\newcommand {\la} {\lambda}
\newcommand{\om}{\omega}
\newcommand{\Om}{\Omega}
\newcommand {\sa} {\sigma}
\newcommand {\Sa} {\Sigma}
\newcommand {\te} {\theta}
\newcommand {\fy} {\varphi}
\newcommand {\ep} {\varepsilon}
\newcommand{\e}{\varepsilon}
\newcommand {\Dl} {\Delta}
\newcommand {\bfep} {{{\bar \varepsilon}}}


\newcommand{\bi}{{{\bf i}}}
\newcommand{\bj}{{{\bf j}}}
\newcommand{\bk}{{{\bf k}}}
\newcommand{\bg}{{{\bf g}}}
\newcommand{\bh}{{{\bf h}}}
\newcommand{\bv}{{{\bf v}}}
\newcommand{\bw}{{{\bf w}}}
\newcommand{\bx}{{{\bf x}}}
\newcommand{\ba}{{{\bf a}}}
\newcommand{\bl}{{{\bf l}}}
\newcommand{\bbw}{{\overline{\bf w}}}
\newcommand{\bal}{{{\bar \alpha}}}
\newcommand{\bbe}{{{\bar \beta}}}

\newcommand {\gbi}{{\ga_\bi}}
\newcommand {\gbj}{{\ga_\bj}}
\newcommand {\gbk}{{\ga_\bk}}


\newcommand{\eS}{{\EuScript S}}
\newcommand{\eC}{{\EuScript C}}
\newcommand{\eP}{{\EuScript P}}
\newcommand{\eT}{{\EuScript T}}
\newcommand{\eG}{{\EuScript G}}
\newcommand{\eK}{{\EuScript K}}
\newcommand{\eF}{{\EuScript F}}
\newcommand{\eZ}{{\EuScript Z}}
\newcommand{\eL}{{\EuScript L}}
\newcommand{\E}{{\EuScript E}}

\newcommand{\hd}{{\hat d}}
\newcommand{\hS}{{\hat S}}
\newcommand{\hG}{{\hat G}}
\newcommand{\hX}{{\hat X}}
\newcommand{\hT}{{\hat T}}

\newcommand{\gS}{{\mathfrak S}}
\newcommand{\gF}{{\mathfrak F}}
\newcommand{\m}{\frak{M}}

\newcommand{\cO}{{\mathcal O}}
\newcommand{\cV}{{\mathcal V}}

\newcommand \bbc {\mathbb{C}}
\newcommand{\bbC}{\mathbb{C}}
\newcommand{\rr}{\mathbb{R}}
\newcommand \nn {\mathbb{N}}
\newcommand \zz {\mathbb{Z}}
\newcommand \rd {\mathbb{R}^d}

\newcommand{\vA}{{\vec {A}}}
\newcommand{\vB}{{\vec {B}}}
\newcommand{\vF}{{\vec {F}}}
\newcommand{\vf}{{\vec {f}}}
\newcommand{\vK}{{\vec {K}}}
\newcommand{\vP}{{\vec {P}}}
\newcommand{\vX}{{\vec {X}}}
\newcommand{\vY}{{\vec {Y}}}
\newcommand{\vx}{{\vec {x}}}
\newcommand{\va}{{\vec {a}}}

\newcommand{\tA}{{\tilde {A}}}
\newcommand{\tB}{{\tilde {B}}}
\newcommand{\tK}{{\tilde {K}}}
\newcommand{\tT}{{\tilde {T}}}
\newcommand{\tU}{{\tilde {U}}}
\newcommand{\tO}{{\tilde {O}}}


\newcommand{\Tbar}{{\overline{T}}}
\newcommand{\Kbar}{{\overline{K}}}
\newcommand{\Mbar}{{\overline{M}}}
\newcommand{\Nbar}{{\overline{N}}}
\newcommand{\Sbar}{{\overline{S}}}
\newcommand{\Xbar}{{\overline{X}}}
\newcommand{\Ybar}{{\overline{Y}}}
\newcommand{\wbar}{{\overline{w}}}
\newcommand{\zbar}{{\overline{z}}}
\newcommand{\gbar}{{\overline{g}}}
\newcommand{\Gbar}{{\overline{G}}}
\newcommand{\Abar}{{\overline{A}}}
\newcommand{\Aw}{{\widetilde{A}}}
\newcommand{\Bw}{{\widetilde{B}}}
\newcommand{\Hw}{{\widetilde{H}}}
\newcommand{\Sw}{{\widetilde{S}}}
\newcommand{\Uw}{{\widetilde{U}}}
\newcommand{\Xw}{{\widetilde{X}}}
\newcommand{\Tw}{{\widetilde{T}}}
\newcommand{\Gw}{{\widetilde{G}}}
\newcommand{\gw}{{\widetilde{g}}}

\newcommand{\shift}{\mathsf{s}}
\newcommand{\sS}{\mathsf{S}}

\def \Div {\mathop{\rm div}\nolimits}
\def \tg {\mathop{\rm tg}\nolimits}
\def \ctg {\mathop{\rm ctg}\nolimits}
\def \arctg {\mathop{\rm arctg}\nolimits}
\def \arcctg {\mathop{\rm arcctg}\nolimits}
\def \arcsin {\mathop{\rm arcsin}\nolimits}
\def \re {\mathop{\rm Re}\nolimits}
\def \im {\mathop{\rm Im}\nolimits}
\def \Arg {\mathop{\rm Arg}\nolimits}
\def \arg {\mathop{\rm arg}\nolimits}
\def \const {\mathop{\rm const}\nolimits}
\def \sh {\mathop{\rm sh}\nolimits}
\def \ch {\mathop{\rm ch}\nolimits}
\def \ind {\mathop{\rm ind}\nolimits}
\def \res {\mathop{\rm Res}\nolimits}
\def \diam {\mathop{\rm diam}\nolimits}
\def \card {\mathop{\rm card}\nolimits}
\def \Id {\mathop{\rm Id}\nolimits}
 \def \fix {\mathop{\rm fix}\nolimits}
 \def \Lip {\mathop{\rm Lip}\nolimits}
 \def \Re {\mathop{\rm Re}\nolimits}
 \def \Im {\mathop{\rm Im}\nolimits}
 \def \max {\mathop{\rm max}\limits}
\def \min {\mathop{\rm min}\limits}
 \newcommand \ls {\mathop{\rm lim\ sup}\limits}
 \newcommand \li {\mathop{\rm lim\ inf}\limits}

\newcommand{\VEC}{\overrightarrow}
\newcommand{\IN}{{\subset}}
\newcommand{\NI}{{\supset}}
\newcommand \dd  {\partial}
\newcommand {\mmm}{{\setminus}}
\newcommand{\probel}{\vspace{.5cm}}
\newcommand{\8}{{\infty}}
\newcommand{\vse}{$\blacksquare$}

 \newtheorem{teo}{\sc Theorem}[section]
 \newtheorem{sled}[teo]{\sc Corollary}
 \newtheorem{lem}[teo]{\sc Lemma}

\newtheorem{opr}[teo]{\sc Definition}

Let $\eS$ be a system $\{S_1,...,S_m\}$ of injective contraction
maps of a complete metric space $(X,d)$ to itself and let $K$ be
it's {\em invariant set},  i.e. such a nonempty compact set $K$
that satisfies $K=\bigcup\limits_{i=1}^m S_i(K)$. The set $K$ is
also called {\em the attractor} of the system $\eS$. A natural
construction allowing to obtain the systems $\eS$ with a connected
(and therefore arcwise connected) invariant set is called a
self-similar zipper and it goes back to the works of Thurston
\cite{Thu} and Astala \cite{Ast} and was analyzed in detail by
Aseev, Kravtchenko and Tetenov in \cite{AKT}. Namely,
\begin{opr}
 A system $\eS=\{S_1,...,S_m\}$ of injective
contraction maps of complete metric space $X$ to itself is called
a zipper with vertices $(z_0,...,z_m)$
 and signature $\vec\ep=(\ep_1,...,\ep_m)\in\{0,1\}^m$   if for
 any $j=1,...,m$ the following equalities hold: 1.
 $S_j(z_0)=z_{j-1+\ep_j}$;  2. $S_j(z_m)=z_{j-\ep_j}$.\end{opr}

If the maps $S_i$ are similarities (or affine maps) the zipper is
called self-similar (correspondingly self-affine).

We shall call the points $z_0$ and $z_m$  {\em the initial} and
{\em the final} point of the zipper respectively.

The simplest example of a self-similar zipper may be obtained if
we take a partition $P$, $0=x_0<x_1<\ldots<x_m=1$ of the segment
$I=[0,1]$ into $m$ pieces and put
$T_i=x_{i-1+\ep_i}(1-t)+x_{i-\ep_i}t $. This zipper $\{T_1,\ldots,
T_m\}$ will be denoted by $\eS_{P,\vec\ep}$.

\begin{teo} {\rm( see \cite{AKT})}. For any zipper $\eS=\{S_1,...,S_m\}$ with
vertices $\{z_0, \ldots, z_m\}$ and signature $\vec\ep$ in a
complete metric space $(X,d)$ and for any partition
$0=x_0<x_1<\ldots<x_m=1$ of the segment $I=[0,1]$ into $m$ pieces
there exists  unique map $\ga:I\to K(\eS)$ such that for each
$i=1,...,m$, $\ga(x_i)=z_i$ and $S_i\cdot\ga=\ga\cdot T_i$ (where
$T_i\in \eS_{P,\vec\ep} $). Moreover, the map $\ga$ is H\"older
continuous.\end{teo}

The mapping $\ga$ in the Theorem is called a {\em linear
parametrization} of the zipper $\eS$. Thus, the attractor $K$ of
any zipper $\eS$ is an arcwise connected set, whereas the linear
parametrization $\ga$ may be viewed as a self-similar Peano curve,
filling the continuum $K$.
 \vskip .5cm {\bf Some Peano curves.}

a) The attractor $K$ of a self-similar zipper $\eS$ with vertices
$(0, 0)$, $ (1/4, \sqrt{3}/4)$, $ (3/4, \sqrt{3}/4)$, $ (1, 0)$
and signature $(1,0,1)$ is  the Sierpinsky gasket.
\begin{figure}[h]
     \epsfig{file=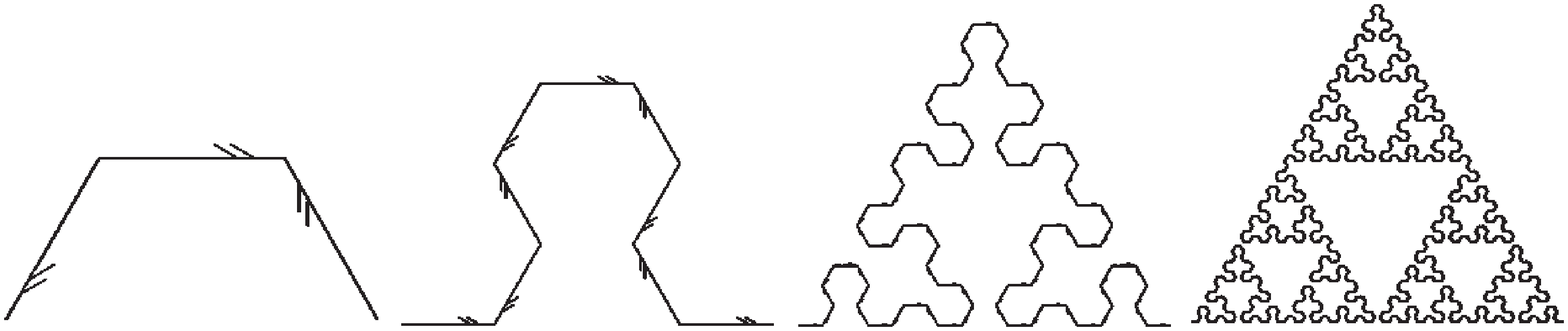,angle=0,scale=.4}
  \caption{\tt  1,2,4,and 8 iterations in the construction of the Peano curve for Sierpinsky gasket.}
  \label{fig:serp-iter}
\end{figure}
\vskip .5cm b) A self-similar zipper with vertices $(0,0)$,
$(0,1/2)$, $(1/2,1/2)$, $(1,1/2)$, $(1,0)$ and signature (1,0,0,1)
produces a self-similar Peano curve for the square $[0,1]\times
[0,1] $
\begin{figure}[h]
\epsfig{file=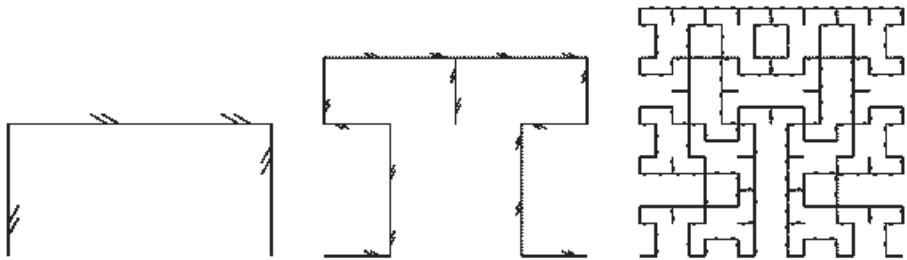,angle=0,scale=.8}\label{f2}
       \caption{\tt  Iterations for square-filling Peano curve.}
  \label{fig:peano-iter}
\end{figure}

\vskip .5cm c) A self-similar zipper with vertices
$(0,0),(0,1/3)$, $(1/3,1/3),(1/3,2/3)$, $(1/3,1),(2/3,1)$,
$(2/3,2/3),(2/3,1/3)$, $(2/3,0),(1,0)$ and signature (0,1,0,
0,1,0, 0,1,0) gives a Peano curve for Sierpinsky carpet.

\vskip .5cm d) The attractor of a zipper with vertices $(0,0),$
$(1,0),$ $(1,1)$, $(1,2)$, $(2,2),$ $(2,1),$ $(2,0),$ $(3,0)$ and
signature (0,0,1,1,1,0,0) is a dendrite.

\begin{figure}[h]
\epsfig{file=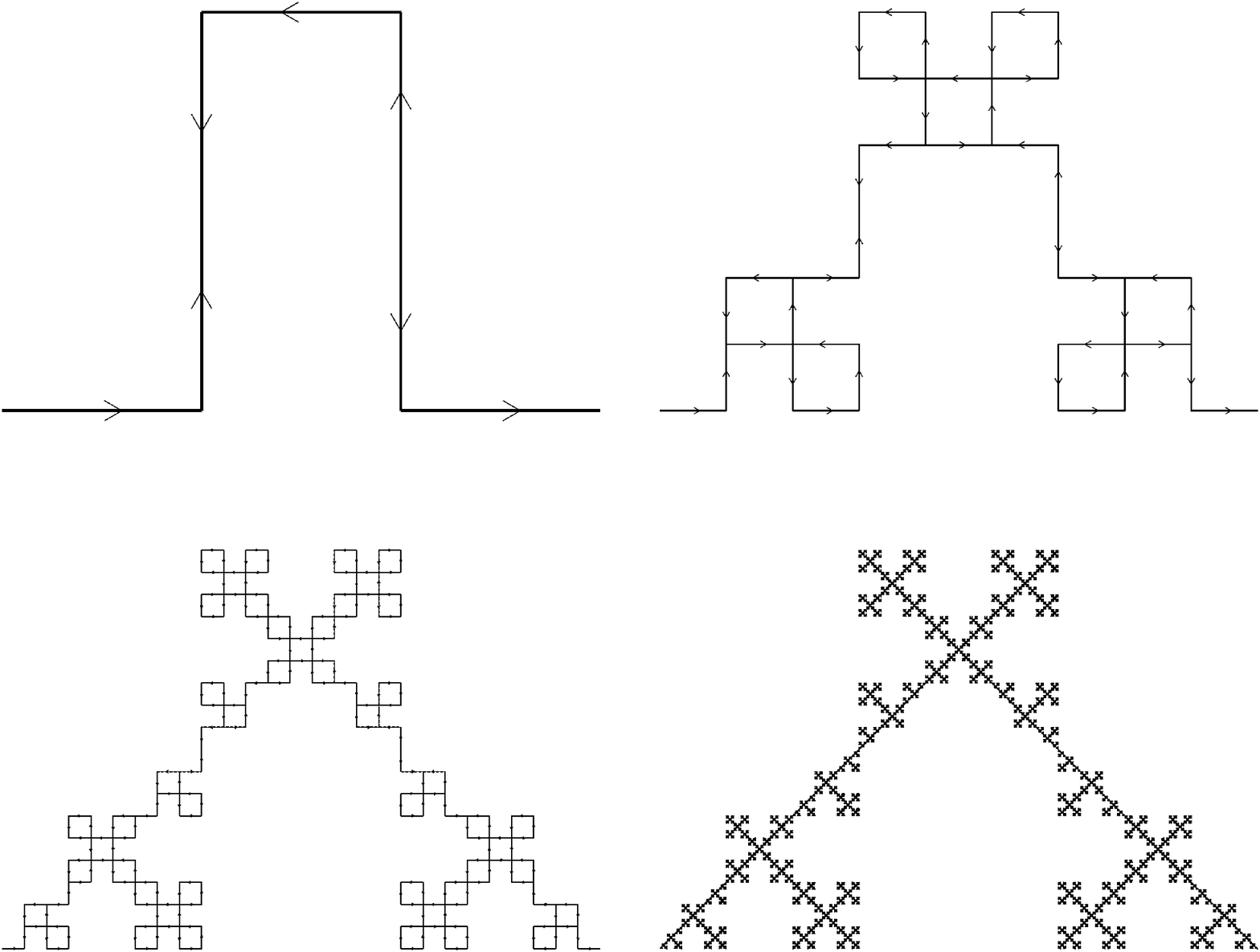,angle=0,scale=.15}\label{f2}
       \caption{\tt  A zipper whose attractor is a dendrite.}
  \label{fig:zipdendr}
\end{figure}

\vskip .5cm{\bf The main example.}

The following example shows that there do exist self-similar
continua which cannot be represented as an attractor of a
self-similar zipper.

Let $\eS$ be a system of contraction similarities $g_k$ in $\rr^2$
where $S_2(\vec x)=\vx/2+(2,0)$, and $S_k(\vec x)=\vx/4+\va_k$
where $\va_k$ run through the set $\{$ $(0,0)$,  $(3,0)$,
$(1,2h)$, $(3/2,3h)\}, h=\sqrt{3}/2$ for $k=1,3,4,5$. Let $K$ be
the attractor of the system $\eS$  and $T$ -- the Hutchinson
operator of the system $\eS$ defined by
$T(A)=\bigcup\limits_{j=1}^5 S_j(A)$.

\begin{figure}[h]
\epsfig{file=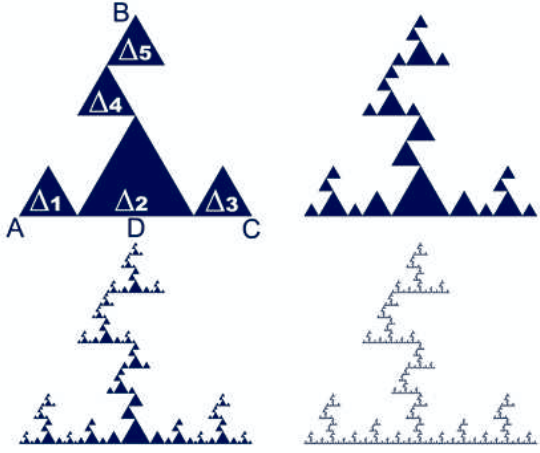,angle=0,scale=.5}\label{f2}
       \caption{\tt  Iterations for the example.}
  \label{fig:peano-iter}
\end{figure}

\vskip .5cm We shall use the following notation: By $\Delta$ we
denote the triangle with vertices $A=(0,0)$, $B=(2,2\sqrt{3})$ and
$C=(4,0)$. The point $(2,0)$ is denoted by $D$. For a multiindex
$\bi=i_1...i_k$ we denote $S_\bi=S_{i_1}...S_{i_k}$,
$\Delta_\bi=S_\bi(\Delta)$, $K_\bi=S_\bi(K)$, $A_\bi=S_\bi(A)$,
etc.

\vskip .5cm 1. {\em The set $K$ is a dendrite.} The way the system
$\eS$ is defined (see \cite[Thm.1.6.2]{Kig}) guarantees the
arcwise connectedness of $K$. Since for each $n$ the set
$T^n(\Delta)$ is simply-connected, the set $K$ contains no cycles
and therefore $K$ is a dendryte. Each point of $K$ has the order 2
or 3. If a point $x$ has the order 3, it is an image $S_\bi(D)$ of
the point $D$ for some multiindex $\bi$. Any path in $K$
connecting a point $\xi\in J$ with a point $\eta\in \Delta_\bi,
\bi=4,5,24,25,224,225,..$, passes through the point $D$.

\vskip .5cm 2. {\em Each non-degenerate line segment $J$ contained
in $K$, is parallel to x axis and is contained in some maximal
segment in $K$ which has the length $4^{1-n}$.}

Consider a non-degenerate linear segment $J\IN K$. There is such
multiindex $\bi$, that $J$ meets the boundary of $S_\bi(\Delta)$
in two different points which lie on different sides of
$S_\bi(\Delta)$ and do not lie in the same subcopy of $K_\bi$.
Then $J'=g_{\bi}^{-1}(J\cap K_\bi)$ is a segment in $K$ with the
endpoints lying on different sides of $D$ which is not contained
in neither of subcopies $K_1,...,K_5$ of $K$. Then $J'=[0,4]$.
Since a part of $J$ is a base of some triangle $S_\bi(\Delta)$,
the length of the maximal segment in $K$ containing $J$ is
$4^{1-n}$ where $n\le |\bi|$.

\vskip .5cm 3. {\em Any injective affine mapping $f$ of $K$ to
itself is one of the similarities $S_\bi=S_{i_1}\cdot...\cdot
S_{i_k}$.} Since $f$ maps $[0,4]$ to some $J\subset S_\bi([0,4])$
for some $\bi$, it is of the form $f(x,y)=(ax+b_1y+c_1,b_2y+c_2)$,
with positive $b_2$. Choosing appropriate composition
$S_\bi^{-1}\cdot f\cdot S_\bj(K)$ we obtain a map of $K$ to itself
sending $[0,4]$ to some subset of $[0.4]$.

Therefore we may suppose that $f(x,y)=(ax+b_1y+c_1,b_2y)$, and
that the image $f(\Delta)$ is contained in $\Delta$ and is not
contained in any $\Delta_i, i=1,...,5$.

If $f(B)\in \Delta_\bi, \bi=4,5,24,25$, then, since every path
from $J$ to $f(B)$ passes through $D$,   $f(D)=D$ and therefore
$c_1=2-a$.

If $f(B)\in \Delta_i, i=4,5$, then $1/2\le b_2\le 1$. In this case
$y-$coordinates of the points $f(B_1),f(B_3)$ are greater than
$\sqrt{3}/4$, so they are contained in $\Delta_1$ and $\Delta_3$,
therefore the map $f$ either keeps the points $D_1, D_3$ invariant,
or transposes them. In each case $|a|=1$ and $f(\{A,C\})=\{A,C\}$.
If in this case $f(B)\neq B$, then $f(A_4)$ cannot be contained in
$T(\Delta)$. The same argument shows that if $f(B)=B$, then
$f(A)\neq C$. Therefore $f=\Id$.

Suppose $f(B)\in \Delta_i, i=24,25$ and $a>1/2$. Then  the points
$f(B_1),f(B_3)$ are contained in $\Delta_1$ and $\Delta_3$,
therefore the map $f$ either keeps the points $D_1, D_3$
invariant, or transposes them, so $|a|=1$ and
$f(\{A,C\})=\{A,C\}$. Considering the intersections of the line
segments $[A,f(B)]$ and $[f(B),C]$ with the boundary of
$T(\Delta)$ and $T^2(\Delta)$ we see that either  $f(A_4)$ or
$f(C_5)$ is not contained in $T^2(\Delta)$, which is impossible.

Therefore, either  $a\le 1/2$ or $f=\Id$ . The first means that
$f(\Delta)\subset \Delta_2$, which contradicts the original
assumption, so $f=\Id$.

\vskip .5cm 4. {\em The set $K$ cannot be an attractor of a
zipper.}
 Let
$\Sigma=\{\varphi_1,...,\varphi_m\}$ be a zipper whose invariant
set is $K$. Let $x_0, x_1$ be the initial and final points of the
zipper $\Sigma$. Let $\gamma$ be a path in $K$ connecting $x_0$
and $x_1$. Since for every $i=1,...,m$ the map $\varphi_i$ is
equal to some $S_\bj$, the sets $\varphi_i(K)$ are the subcopies
of $K$, therefore for each $i$ at least one the images
$\varphi_i(x_0),\varphi_i(x_1)$ is contained in the intersection
of $\varphi_i(K)$ with adjacent copies of $K$. Consider the path
$\tilde\gamma=T_\Sigma(\gamma)=\bigcup\limits_{i=1}^m
\varphi_i(\gamma)$. It starts from the point $x_0$, ends at $x_1$
and passes through all copies $K_j$ of $K$. Each of the points
$C_1=A_2$, $C_2=A_3$, $B_2=C_4$ and $B_4=A_5$ splits $K$ to two
components, therefore is contained in $\tilde\gamma$ and is a
common point for the copies  $\varphi_i(\gamma),$ $
\varphi_{i+1}(\gamma)$ for some $i$. Therefore one of the points
$x_0,x_1$ must be $A$,
 one of the points $x_0,x_1$ must be $B$, and one of the points $x_0,x_1$ must be
 $C$, which is impossible.

 \end{document}